\documentclass[runningheads]{llncs}
\usepackage[T1]{fontenc}
\usepackage{newtxtext}
\usepackage[varvw]{newtxmath}
\usepackage{graphicx}
\usepackage{amsmath}
\usepackage{algpseudocode}
\usepackage{algorithm}
\usepackage{color}
\usepackage[numbers]{natbib}
\usepackage{hyperref}

\begin{document}
\title{To Race or Not to Race? Exploiting Performance Variability in Mixed-Integer Programming via Root-Gated Resource Allocation}
\titlerunning{To Race or Not to Race?}


\author{Mourchid Adegbindin \and Brenda Dietrich, \and Jamol Pender}
\authorrunning{Adegbindin et al.}
\institute{Department of Operations Research and Information Engineering\\
Cornell University, Ithaca, NY, USA, 14853\\ \email{ma757@cornell.edu, bld34@cornell.edu, jjp274@cornell.edu}
}

\maketitle

\begin{abstract}
Mixed-integer programming (MIP) solver runtimes vary
substantially across random seeds. Parallel replication
exploits this variability but can waste computational
resources when applied indiscriminately. We investigate
whether information from an ordinary solve can guide
this allocation decision.

We formulate replication under a finite deadline,
establish an exact relationship between runtime dispersion
and the value of independent searches, and analyze delayed
replication and imperfect allocation decisions. Using
101,615 seeded SCIP runs, we examine the association
between root-LP dual degeneracy and replication value.
We then introduce a non-invasive root-gated policy that
preserves the original solver trajectory and launches
three additional seeded processes when a root
dual-degeneracy indicator reaches a fixed threshold.

In a prospective two-hour experiment on 130 paired
instances, root-gated racing proves optimality on 77
instances (59.23\%), versus 73 (56.15\%) for both Single and unconditional Race-4,
while reducing active worker-process occupancy by
43.00\% relative to unconditional racing. These findings
demonstrate the promise of selective parallel resource
allocation without solver-source modifications.

\end{abstract}

\section{Introduction}
\label{sec:introduction}
Mixed-integer programming (MIP) is an essential tool for
planning, scheduling, and resource allocation. In many
applications, an optimization run must deliver a solution
within a prescribed time window under a limited computational
budget. The usefulness of the resulting decision depends
directly on the reliability of that computation.

Yet modern MIP solvers exhibit substantial runtime variability
across random seeds: identical models, settings, and
hardware can yield markedly different branch-and-bound
trajectories, with one run proving optimality while another
reaches the time limit. This unpredictability complicates
resource planning: insufficient capacity risks missed
deadlines, while excessive allocation wastes computation.

Independent seeded searches can improve the chance of a
timely proof, but consume additional processor time,
memory, and computing capacity, often without benefit
on less seed-sensitive instances. The central question is
therefore: \emph{When does information available during an
ongoing MIP solve justify allocating additional processors
to alternative search trajectories?} We address this
decision under a finite deadline, accounting explicitly
for wall-clock time and aggregate computational effort.

\subsection{Approach and Related Work}
\label{sec:intro_prior_work}

Performance variability in MIP is well documented
\cite{danna2008,KochEtAl2011_MIPLIB2010,
LodiTramontani2013,halSortModel}.
Adegbindin et al.~\cite{Adegbindin2026VariabilityScore}
quantify this variability across randomized SCIP runs and
examine its relationship with instance structure and early
solver information. Their work establishes that performance
variability is not merely experimental noise, but a
measurable feature of MIP computation with potential
operational consequences.

The exploitation of variability through randomized search
is well established. Algorithm portfolios exploit
performance differences across alternative strategies
\cite{gomes2001portfolios}. Fischetti and Monaci
\cite{fischetti2014erraticism} develop a bet-and-run
approach that conducts short randomized searches and
selects a promising trajectory for continuation.
Fischetti et al.~\cite{fischetti2016random} investigate
randomized sampling in branch-and-cut, including
information collected from diversified root searches.
Carvajal et al.~\cite{carvajal2014parallel} study parallel
branch-and-bound trees that exchange information during
execution. Commercial solvers also support concurrent
MIP optimization through independent searches
\cite{gurobiConcurrentMIP}.

Other methods use information collected during execution
to intervene in the search itself. Anderson et
al.~\cite{anderson2019clairvoyant} use online estimates
of remaining branch-and-bound tree size to guide restart
decisions. Adegbindin et
al.~\cite{Adegbindin2026LPStateSteering} develop an
LP-state steering framework that modifies the state
inherited by branch-and-bound. These approaches
demonstrate the value of intervening in an ongoing
computation, but address different control decisions.

Our work addresses a complementary question. Rather
than modifying the search trajectory, we ask whether
information produced by an ordinary solve can determine
when additional search trajectories are worth their
computational cost. One solver process begins normally,
and the controller observes its first solved root linear
program (LP). At this checkpoint, the controller decides
whether to continue with the original process alone or
allocate additional processors to independent seeded
searches. The original trajectory is preserved in either
case.

The decision relies on a simple dual-degeneracy statistic:
the fraction of eligible nonbasic LP columns with zero or
near-zero reduced cost. Gamrath, Berthold, and Salvagnin
\cite{gamrath2020degeneracy} study this classical
degeneracy measure; we investigate its value as an early
resource-allocation signal. A fixed threshold determines
whether three additional seeded processes launch.
The controller operates through SCIP's existing interface
and standard operating-system process controls, without
modifying solver source code or exchanging information
between replicas.

This design creates a tradeoff absent from immediate
replication. Additional processes start only after the
root LP becomes available, and their computational cost
must be justified against the remaining deadline.
We analyze this tradeoff through the distributional
value of independent replication, delayed launches,
and imperfect allocation decisions.

\noindent\textbf{Contributions.}
We make three contributions. First, we develop a
finite-deadline analysis of randomized replication,
with explicit accounting for delayed launches,
processor consumption, and imperfect allocation
decisions. Second, using 101,615 seeded SCIP runs,
we examine whether root-LP dual degeneracy identifies
instances for which additional search trajectories
have value, and quantify the opportunities captured
and missed by a fixed threshold. Third, we implement
a solver-external controller that preserves the
original trajectory while selectively launching
independent seeded searches. On 130 prospectively
paired instances at a two-hour deadline, Root-Gated
proves optimality on 77 instances, against 73 for
both single-seed solving and unconditional Race-4,
while reducing active worker occupancy by 43.00\%
relative to Race-4.

\section{When Does Racing Pay?}
\label{sec:racing_theory}

\subsection{Replication and Computational Cost}
\label{sec:replication_model}

Consider a fixed MIP instance, solver configuration, and
computing environment. Let $X\geq 0$ denote the time required
to prove optimality under a randomly selected seed, with
$X=\infty$ if the trajectory never proves optimality. In an
idealized replication model, $X_1,\ldots,X_m$ are independent
copies of $X$.

When all $m$ processes start simultaneously and terminate at
the first proof, the completion time is
$M_m=\min_{1\leq i\leq m}X_i$. For a common deadline $T$, let
$M_m^T=\min\{M_m,T\}$ and $\bar F(t)=\Pr(X>t)$. Independence
yields
\begin{equation}
    \mathbb E[M_m^T]
    =\int_0^T\bar F(t)^m\,dt,
    \qquad
    \mathbb E[C_m^T]
    =m\int_0^T\bar F(t)^m\,dt,
    \label{eq:replication_survival}
\end{equation}
where $C_m^T=mM_m^T$ denotes aggregate processor-slot
consumption under one continuously active processor per
replica.

For two replicas, suppose $0<\mu=\mathbb E[X]<\infty$.
Define the runtime Gini coefficient by
\[
    G_X=\frac{\mathbb E|X_1-X_2|}{2\mu}.
\]
The identity
\[
    \min\{X_1,X_2\}
    =\frac{X_1+X_2-|X_1-X_2|}{2}
\]
gives
\begin{equation}
    \mathbb E[M_2]=\mu(1-G_X),
    \qquad
    \frac{\mathbb E[C_2^\infty]}{\mu}
    =2(1-G_X).
    \label{eq:gini_replication}
\end{equation}
Consequently,
\begin{equation}
    \mathbb E[C_2^\infty]<\mathbb E[X]
    \quad\Longleftrightarrow\quad
    G_X>\frac12.
    \label{eq:gini_condition}
\end{equation}

Thus, sufficiently dispersed runtimes can make replication
less expensive in expected processor time as well as faster
in expected wall-clock time. Variance alone does not
determine this benefit. Under a finite deadline, the same
identities apply to $X^T=\min\{X,T\}$ through its mean and
Gini coefficient.

\subsection{Delayed Replication under a Finite Deadline}
\label{sec:delayed_replication}

A baseline solver process starts at time zero and reaches
its first solved root LP at elapsed time $a<T$. The
remaining wall-clock budget is $b=T-a$. The analysis is
conditional on reaching this checkpoint before the
deadline. If no valid root observation becomes available,
the implemented policy continues the baseline process
without diversification.

At the checkpoint, the controller chooses between
continuing the baseline alone and launching $m-1$
additional seeded processes, where $m\geq 2$. Let $R$ denote the remaining
time to proof for the baseline process in isolation. Let
$X_1,\ldots,X_{m-1}$ denote the completion times of the
alternative processes, measured from their effective
start times, and put
$Y=\min_{1\leq i\leq m-1}X_i$.

The additional processes become active after a common,
fixed delay $h$, where $0\leq h<b$. Let $H\geq 0$ denote
additional processor-equivalent work incurred during
launch, monitoring, and termination that is not included
in the active-replica consumption below. We assume
$\mathbb E[H]<\infty$.

Assuming that diversification does not alter the baseline
trajectory, the remaining elapsed times under continuation
and diversification are
\begin{equation}
    W_0=\min\{R,b\},
    \qquad
    W_1=\min\{R,h+Y,b\}.
    \label{eq:delayed_times}
\end{equation}

We evaluate both actions using processor-slot consumption,
a wall-clock weight $\lambda\geq 0$, and a penalty $P>0$
for failure to prove optimality by the deadline. The
remaining cost of continuation is
\begin{equation}
    J_0=(1+\lambda)W_0
        +P\mathbf 1_{\{R>b\}}.
    \label{eq:continue_cost}
\end{equation}
For diversification, the additional replicas consume
processor time after delay $h$:
\begin{equation}
\begin{aligned}
    J_1={}&(1+\lambda)W_1
       +(m-1)(W_1-h)_+ +H\\
       &+P\mathbf 1_{\{R>b,\ h+Y>b\}},
\end{aligned}
    \label{eq:diversify_cost}
\end{equation}
where $(x)_+=\max\{x,0\}$. The realized incremental cost
of diversification is therefore
\begin{equation}
\boxed{
\begin{aligned}
    \Delta:=J_1-J_0
    ={}&(1+\lambda)(W_1-W_0)\\
       &+(m-1)(W_1-h)_+ +H\\
       &-P\mathbf 1_{\{R>b,\ h+Y\leq b\}}.
\end{aligned}
}
\label{eq:incremental_cost}
\end{equation}

Let $\mathcal I$ denote the information available at the
checkpoint, including $a$ and the fixed model parameter
$h$. Define
\[
    S_0(t)=\Pr(R>t\mid\mathcal I),
\]
and
\[
    S_1(t)
    =\Pr(R>t,\ h+Y>t\mid\mathcal I).
\]
For $t<h$, we have $S_1(t)=S_0(t)$.

\begin{proposition}[Conditional value of delayed replication]
\label{prop:delayed_value}
Under the preceding assumptions,
\begin{equation}
\boxed{
\begin{aligned}
\mathbb E[\Delta\mid\mathcal I]
={}&\mathbb E[H\mid\mathcal I]\\
&+(m-1)\int_h^b S_1(t)\,dt\\
&-(1+\lambda)\int_h^b
    [S_0(t)-S_1(t)]\,dt\\
&-P[S_0(b)-S_1(b)].
\end{aligned}
}
\label{eq:conditional_replication_value}
\end{equation}
\end{proposition}

The proof is given in Appendix~\ref{sec:appendix_proofs}.
Diversification has lower conditional expected cost
precisely when
$\mathbb E[\Delta\mid\mathcal I]<0$.

The model assumes simultaneous activation of the
additional replicas and continuous activity until the
race terminates. The implementation permits staggered
starts and individual worker termination; its resource
consumption is measured from observed process lifetimes.
Moreover, active worker occupancy differs from capacity
reserved by a fixed-size batch allocation.

\subsection{Allocation Decisions with an Imperfect Signal}
\label{sec:imperfect_gate}

Let $G\in\{0,1\}$ be an $\mathcal I$-measurable decision
rule, where $G=1$ indicates diversification. Its realized
cost is $J_G=J_0+G\Delta$. Consequently,
\begin{equation}
\begin{aligned}
    \mathbb E[J_G]<\mathbb E[J_0]
    &\iff \mathbb E[G\Delta]<0,\\
    \mathbb E[J_G]<\mathbb E[J_1]
    &\iff \mathbb E[(1-G)\Delta]>0.
\end{aligned}
\label{eq:gate_comparisons}
\end{equation}

Let
$\delta(\mathcal I)=\mathbb E[\Delta\mid\mathcal I]$
and define the optimal checkpoint rule
\[
    G^\star
    =\mathbf 1_{\{\delta(\mathcal I)<0\}}.
\]

\begin{proposition}[Cost of an imperfect allocation signal]
\label{prop:imperfect_cost}
For every $\mathcal I$-measurable binary gate $G$,
\begin{equation}
\boxed{
\mathbb E[J_G]-\mathbb E[J_{G^\star}]
=
\mathbb E\left[
|\delta(\mathcal I)|
\mathbf 1_{\{G\ne G^\star\}}
\right].
}
\label{eq:gate_regret}
\end{equation}
\end{proposition}

The proof is given in Appendix~\ref{sec:appendix_proofs}.
The regret from an allocation error depends on its
expected cost, rather than its frequency alone.

For comparison, let $B_q$ be a binary random decision
independent of $(J_0,J_1,G)$, with
$\Pr(B_q=1)=q=\mathbb E[G]$. Then
\begin{equation}
    \mathbb E[J_G]-\mathbb E[J_{B_q}]
    =\operatorname{Cov}(G,\Delta).
    \label{eq:gate_random_comparison}
\end{equation}

The optimal rule $G^\star$ depends on conditional runtime
and cost distributions unavailable for a new instance.
We therefore restrict attention to an inexpensive
structural statistic obtained during the baseline solve.
The resulting threshold policy approximates the
conditional cost decision; its effectiveness depends
on the association between the statistic and the
marginal value of replication. Section~
\ref{sec:stagewise_variability} examines this association
using historical solver runs.

\subsection{The Time-Zero Racing Benchmark}
\label{sec:timezero_benchmark}

Let $X_0=a+R$ denote the baseline trajectory's total
completion time, conditional on reaching the checkpoint,
and let $Y$ denote the alternative portfolio's completion
time under the same isolated seeded trajectories.
Under a coupling that preserves these trajectories
across policies, immediate time-zero racing has
completion time
\[
    D_{\mathrm{race0}}=\min\{X_0,Y\}.
\]
For a root-gated policy,
\begin{equation}
    D_G=
    \begin{cases}
      X_0, & G=0,\\
      \min\{X_0,a+h+Y\}, & G=1.
    \end{cases}
\end{equation}
Since $a+h\geq 0$,
\begin{equation}
    D_{\mathrm{race0}}\leq D_G
    \qquad\text{almost surely}.
    \label{eq:timezero_boundary}
\end{equation}
Thus, under this coupling, immediate racing provides a
lower bound on completion time and an upper bound on
deadline proof probability for root-gated racing. This
ordering need not hold when concurrent execution alters
the individual solver trajectories through resource
contention or other interactions.

The resource-allocation objective is to recover useful
deadline protection while reducing unnecessary parallel
computation. The implemented decision rule is
\begin{equation}
    G=\mathbf 1_{\{d_{\mathrm{root}}\geq d^\star\}},
    \label{eq:implemented_gate}
\end{equation}
when a valid root statistic is available, and $G=0$
otherwise. Here $d_{\mathrm{root}}$ is the root
dual-degeneracy indicator and $d^\star$ is a fixed
threshold. Its effectiveness is evaluated empirically,
rather than assumed from the replication theory.


\section{When Does Variability Become Actionable?}
\label{sec:stagewise_variability}

\subsection{The Seeded-Run Corpus}
\label{sec:log_corpus}

We examine 101,615 seeded SCIP runs covering
1,017 distinct instances from the MIPLIB 2017
collection \cite{GleixnerEtAl2021_MIPLIB2017},
under a nominal four-hour time limit. Of these, 1,016 instances have 100
recorded runs and one has 15. The records contain
terminal statuses, solving times, random seeds,
branch-and-bound statistics, and early solver
measurements.

Among these runs, 42,877 terminate optimally,
56,452 reach the time limit, 1,704 have other
statuses, and 582 lack parsed status. First-LP
basis and reduced-cost counts are available for
92,548 runs, with a defined root-degeneracy ratio
in 92,048. Appendix~\ref{sec:empirical_appendix}
provides the full audit.

We count an optimality proof only when the
recorded status is optimal and its completion
time does not exceed deadline $T$; time-limited
runs are right-censored. The Weighted Variability
Score (WVS) of \cite{Adegbindin2026VariabilityScore}
describes cross-seed variability but is not
an input to the online allocation policy.

\subsection{Empirical Value of Replication}
\label{sec:replication_opportunity}

We first isolate the value of replication
independently of any root-based decision rule.
For each instance, we form disjoint blocks of
four recorded seeds: one baseline and three
alternative trajectories. An ideal four-seed
race proves optimality by deadline $T$ whenever
at least one seed does. Its completion time
is the minimum of the four trajectories,
capped at $T$. This time-zero counterfactual
excludes launch delays, overhead, and hardware
contention.

The audit yields 24,798 eligible blocks from
994 instances. At four hours, ideal Race-4
increases proof probability from 43.04\%
to 45.67\%, recovering 654 additional proofs.
Among 113 retrospectively identified mixed-outcome
instances, the 2,818 eligible blocks yield
47.48\% and 70.69\%, respectively: a 23.21
percentage-point gain. At 15 minutes, the
corresponding all-block fractions are 31.08\%
and 33.51\%.

This reliability gain carries substantial
cost: mean ideal processor consumption rises
from 2.44 to 9.35 processor-hours per block.
The historical proof curves appear in
Appendix~\ref{sec:empirical_appendix}.
Thus, replication creates valuable proof
opportunities, but unconditional deployment
expends substantial resources on instances
with little additional benefit.

\subsection{Root Dual-Degeneracy Indicator}
\label{sec:root_degeneracy}

We investigate whether the first solved root
LP identifies valuable replication opportunities.
Our signal is the classical dual-degeneracy rate
\cite{gamrath2020degeneracy}. Let $\mathcal N$
denote eligible nonbasic LP columns and
$\bar c_j$ their reduced costs. For numerical
tolerance $\varepsilon_{\mathrm{rc}}$, define
\begin{equation}
    d_{\mathrm{root}}
    =
    \frac{
    \left|\left\{
      j\in\mathcal N:
      |\bar c_j|\leq\varepsilon_{\mathrm{rc}}
    \right\}\right|
    }{
    |\mathcal N|
    }.
    \label{eq:droot}
\end{equation}

The historical records yield the proxy
\begin{equation}
    \widehat d_{\mathrm{root}}
    =
    \frac{
      N_{\mathrm{nonbasic,zeroRC}}
    }{
      N_{\mathrm{columns}}-N_{\mathrm{basic}}
    }.
    \label{eq:historical_droot}
\end{equation}

Records with missing counts or nonpositive
denominators remain undefined, never zero.
Both measurements use nonbasic statuses lower,
upper, or zero and reduced-cost tolerance
$10^{-6}$. The historical proxy reconstructs
logged basis counts, whereas the prospective
controller examines eligible LP columns directly;
exact equivalence is unverified. The threshold
$d^\star=0.25$ was fixed heuristically, without
calibration on the historical corpus.

\subsection{Does Root Degeneracy Identify Valuable Races?}
\label{sec:droot_replication_value}

We associate each baseline root statistic
with the subsequent outcomes of its four
seeded trajectories. Let $Q_{ib}^{(1)}$
and $Q_{ib}^{(4)}$ indicate optimality proofs
by deadline for the baseline and ideal race.
The incremental proof benefit is
\begin{equation}
    B_{ib}
       =Q_{ib}^{(4)}-Q_{ib}^{(1)}.
    \label{eq:observed_replication_benefit}
\end{equation}
Here $B_{ib}=1$ identifies a proof obtained
by an alternative seed but missed by the
baseline.

Of the 24,798 eligible blocks, 22,963 from 930
instances have a defined baseline root ratio.
The fixed threshold $d^\star=0.25$ activates
on 11,574 of these blocks (50.40\%) and
identifies 459 of the 654 additional-proof
opportunities (70.18\%). The remaining 195
opportunities lie below threshold; none occur
among blocks with unavailable root measurements.

On the root-observed subset, Single proves
44.22\% of blocks and unconditional ideal
Race-4 proves 47.07\%. Selective time-zero
racing achieves 46.22\% with 5.77 processor-hours
per block, versus 9.18 for Race-4: a 37.15\%
reduction. Independent random activation at
the same 50.40\% rate yields an expected proof
fraction of 45.65\% and 5.82 processor-hours.

\begin{figure}[t]
\centering
\includegraphics[width=\linewidth]
{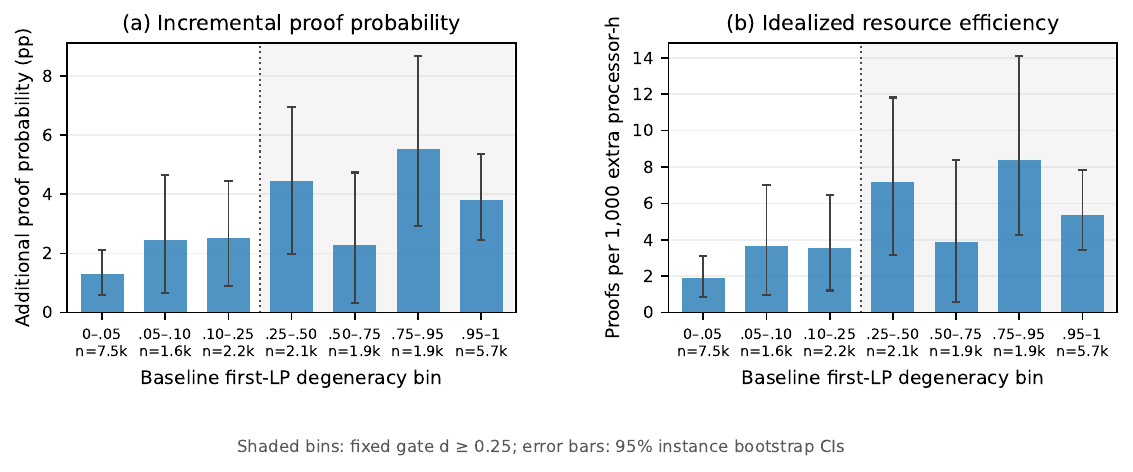}
\caption{Historical first-LP degeneracy and ideal
four-seed replication value at $T=4$ hours:
incremental proof probability (left) and
additional proofs per 1,000 extra processor-hours
(right). Error bars show 95\% instance-bootstrap
intervals; shaded bins satisfy
$d_{\mathrm{root}}\geq0.25$.}
\label{fig:root_replication_value}
\end{figure}

Figure~\ref{fig:root_replication_value}
shows that the association is nonmonotone:
low root degeneracy does not preclude seed
sensitivity, nor does high degeneracy guarantee
beneficial replication. Nevertheless, the
fixed threshold concentrates 70.18\% of
observed proof opportunities within approximately
half the eligible blocks.

These historical comparisons establish the
potential value of selective replication,
not the performance of its delayed implementation.
Section~\ref{sec:computational_evaluation}
evaluates the executable policy prospectively.


\section{Root-Gated Resource Allocation}
\label{sec:root_policy}

\subsection{Root Diversification Gate}
\label{sec:root_gate}

We implement the allocation decision of
Section~\ref{sec:racing_theory} using the root
dual-degeneracy indicator of
Section~\ref{sec:root_degeneracy}. One SCIP process
begins normally; after its first solved root LP,
the controller makes a single allocation decision.
For a fixed threshold $d^\star=0.25$, the rule is
\begin{equation}
G(d_{\mathrm{root}})=
\begin{cases}
0, & d_{\mathrm{root}}<d^\star,\\
1, & d_{\mathrm{root}}\geq d^\star,
\end{cases}
\label{eq:root_allocation_rule}
\end{equation}
where $G=0$ retains the original process alone
and $G=1$ launches three additional seeded searches.
Invalid or unavailable root observations, or a
completed baseline solve, trigger no replicas.
The original trajectory continues undisturbed
in either case. The threshold was fixed
heuristically before the prospective evaluation,
not calibrated on a separate discovery set.

\subsection{Solver-External Implementation and Asynchronous Launch}
\label{sec:external_implementation}

A PySCIPOpt handler captures the
\texttt{FIRSTLPSOLVED} event at the root node.
For an optimal root LP, it computes
$d_{\mathrm{root}}$ from eligible nonbasic columns
with reduced-cost magnitude at most $10^{-6}$.
An atomic JSON record exposes the statistic,
its counts, and validity to an external Python
controller. Neither SCIP source code nor the
active search algorithms are modified.

The maximum portfolio size is $m=4$, with
seeds assigned by
\begin{equation}
    s_{bk}=s_{\mathrm{base}}+4b+k,
    \qquad k=0,1,2,3.
    \label{eq:portfolio_seed_schedule}
\end{equation}
All comparison policies use the same seed block.
Unlike unconditional Race-4, which starts all four
processes immediately, Root-Gated launches replicas
asynchronously after the root observation.
Each receives only the remaining global time
budget; replicas exchange no search information.
The complete controller appears in
Algorithm~\ref{alg:gated_racing} of
Appendix~\ref{sec:empirical_appendix}.

\subsection{Termination and Outcome Accounting}
\label{sec:policy_termination}

The controller evaluates terminal certificates against
the common global deadline using recorded monotonic
completion timestamps. Optimality proofs constitute
the primary endpoint; infeasibility certificates are
recorded separately, and missing terminal records
never count as proofs. Upon recognizing a certificate,
the controller terminates the remaining workers.
Otherwise, processes may continue until their
individual solver time limits expire. Thus, the proof
deadline governs outcome classification, whereas
resource accounting reflects the observed lifetime
of each worker.

For controller-recorded worker launch and stopping
times $s_k,e_k$, aggregate worker-process occupancy
under policy $\pi$ is
\begin{equation}
    C_{\pi}^{\mathrm{slot}}
    =
    \sum_{k\in\mathcal L_\pi}(e_k-s_k),
    \label{eq:portfolio_resource_accounting}
\end{equation}
where $\mathcal L_\pi$ is the set of launched workers.
This measure captures aggregate process occupancy
rather than CPU utilization or scheduler-reserved
capacity.


\section{Computational Evaluation}
\label{sec:computational_evaluation}

We evaluate whether root information supports selective
parallel execution under a common deadline. Unlike the
historical time-zero comparisons, these experiments
execute the three policies prospectively on the cluster.

\subsection{Experimental Design}
\label{sec:discovery_validation}

We use 130 preselected MIPLIB 2017 instances, with one paired
four-seed block per instance (1001-1004). All 390
portfolio trials form complete, deadline-matched
triples at $T=7200$ seconds. The threshold
$d^\star=0.25$ was fixed heuristically before
execution, without calibration on the historical
corpus. An earlier 50-instance run served only
as an engineering shakedown.

Experiments ran through PySCIPOpt 5.4.1 on Seneca's
\texttt{cac\_cpu} partition. Each SCIP worker was
intended to use one core; parallel portfolios reserved
four Slurm CPUs, even when the gate remained closed.
All policies used the same seed block, but their
separate executions could experience different
runtime conditions. We assess uncertainty through
12,000 paired bootstrap resamples of the 130 instances.

\subsection{Policies and Performance Measures}
\label{sec:policy_baselines}

Single executes the baseline seed; Race-4 launches
all four seeds immediately; Root-Gated launches
three additional seeds only when a valid root
measurement reaches $d^\star$. All processes share
the original deadline.

Let $Q_i^\pi$ indicate an optimality proof by the
deadline for policy $\pi$ on instance $i$. We
estimate proof probability and mean worker occupancy
by
\begin{equation}
\widehat p_\pi=\frac{1}{n}\sum_{i=1}^{n}Q_i^\pi,
\qquad
\widehat c_\pi=\frac{1}{n}\sum_{i=1}^{n}
C_{\pi,i}^{\mathrm{slot}},
\qquad n=130.
\label{eq:validation_proof_rate}
\end{equation}
The relative occupancy saving is
\begin{equation}
\mathrm{CostSaving}
=1-\frac{\widehat c_{\mathrm{Gate}}}
{\widehat c_{\mathrm{Race4}}}.
\label{eq:validation_cost_saving}
\end{equation}
Worker occupancy sums active process lifetimes;
portfolio wall time includes controller cleanup.
Neither quantity measures actual CPU utilization
or Slurm-reserved capacity.

\subsection{Main Computational Results}
\label{sec:headline_results}

\begin{table}[t]
\centering
\caption{Prospective results on 130 complete paired
instances at $T=7200$ s. Wall time includes cleanup;
worker-slot seconds measure active process occupancy.}
\label{tab:headline_validation}
\begin{tabular}{lrrr}
\hline
Policy & Proofs (\%) & Wall (s) & Worker-slot (s) \\
\hline
Single & 56.15 & 3588.60 & 3588.59 \\
Race-4 & 56.15 & 3724.71 & 14891.91 \\
Root-Gated & 59.23 & 3499.32 & 8487.73 \\
\hline
\end{tabular}
\end{table}

Root-Gated proves optimality on 77 of 130 instances,
compared with 73 for both Single and Race-4.
Its paired proof advantage is 3.08 percentage
points over either comparator. The gain rests on four
favorable paired discordances and no unfavorable
ones; the exact two-sided McNemar test yields
$p=0.125$. Thus, the observed improvement is
promising but does not establish statistical
superiority under this exact test.

The resource difference is substantial.
Root-Gated reduces mean active worker occupancy
from 14891.91 to 8487.73 seconds per instance:
a 43.00\% saving (95\% bootstrap interval:
34.51-51.39\%), or 6404.18 seconds per trial.
Mean wall time also falls by 225.39 seconds
relative to Race-4 (95\% interval:
117.94-350.37 seconds saved).

\begin{figure}[t]
\centering
\includegraphics[width=\linewidth]{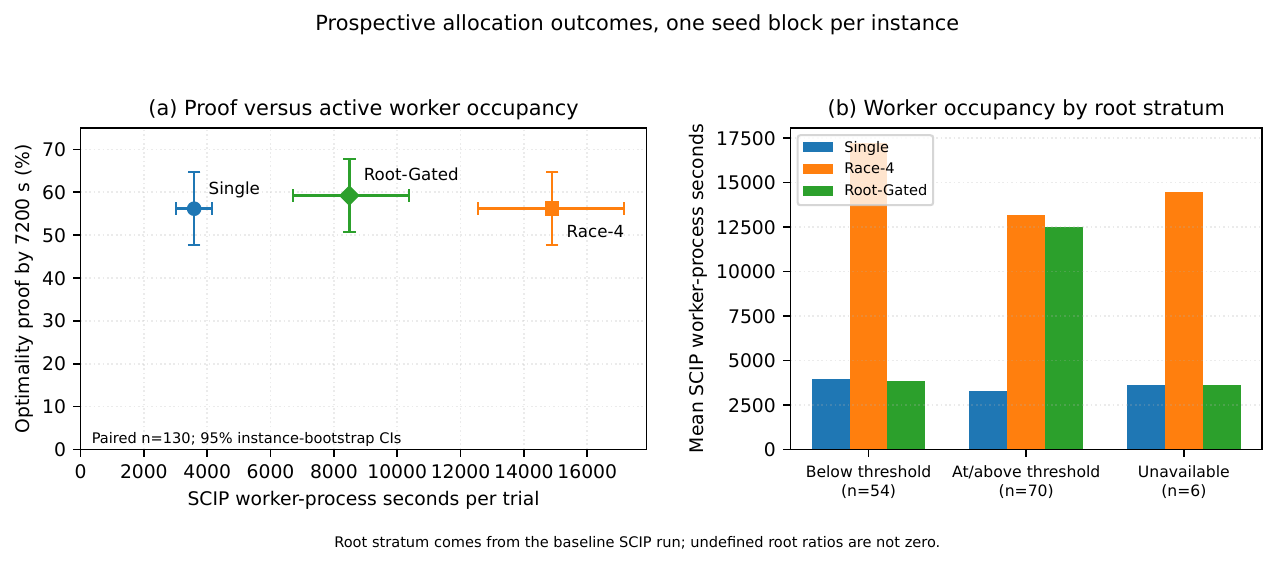}
\caption{Prospective two-hour proof probability
versus active worker occupancy (left), and
worker occupancy by baseline root stratum (right).
Intervals use paired instance bootstraps;
undefined root ratios are not assigned zero.}
\label{fig:performance_resource}
\end{figure}

Active occupancy savings do not translate directly
into reserved cluster capacity. Both parallel
policies reserve four Slurm CPUs throughout
their portfolio jobs. Multiplying mean wall time
by four gives reservation proxies of 14898.84
and 13997.28 CPU-slot seconds for Race-4 and
Root-Gated, respectively: only 6.05\% savings.
Realizing the larger active-occupancy benefit
requires dynamic reallocation of otherwise idle cores.

\subsection{The Value of the Root Gate}
\label{sec:gate_results}

Root-Gated launches four workers on 70 of 130
instances (53.85\%) and retains one on 60.
Of 124 defined root ratios, 70 meet the threshold
and 54 fall below it; six are unavailable.
All gate decisions follow this classification.

For the 54 below-threshold instances, mean worker
occupancy is 3844.81 seconds versus 17200.91
for Race-4. Above threshold, the respective means
are 12487.61 and 13147.61 seconds. The gate thus
concentrates resource savings among instances
without diversification, although these comparisons
do not establish the root statistic's predictive value.

\subsection{Heterogeneity and Failure Cases}
\label{sec:family_results}

Five instances have discordant proof outcomes;
the other 125 agree across policies. Root-Gated
proves four instances missed by each comparator,
but two of its additional proofs relative to
Single and three relative to Race-4 occur with
the gate closed.

The gate decisions and paired proof-outcome
patterns are detailed in
Figure~\ref{fig:gate_decision_audit}
of Appendix~\ref{sec:empirical_appendix}.

At \texttt{fhnw-binpack4-48.mps}, Root-Gated
activates and succeeds where both comparators fail.
At \texttt{k1mushroomi.mps}, it alone succeeds
despite a closed gate. These outcomes expose
execution variability rather than an isolated
causal effect of the threshold. A prospective
random-gate comparator was not executed;
repeated seed blocks and resource-matched controls
remain necessary to assess the signal's
predictive contribution.


\section{Discussion and Practical Implications}
\label{sec:discussion}

Performance variability is usually treated as a source of
computational uncertainty to be measured or mitigated.
Our results suggest a complementary perspective:
variability creates opportunities for selective resource
allocation. Rather than commit four processors before
search begins, a solver can use information from its
first solved root LP to decide whether additional
trajectories warrant their cost.

The prospective experiment illustrates the tradeoff.
Root-Gated obtains 77 optimality proofs against 73
for both comparators, with 43.00\% less active worker
occupancy than unconditional racing. These outcomes
do not, however, establish the causal value of root
dual degeneracy: several discordant proofs occur
without gate activation. The historical analysis
motivates the decision rule, while its predictive
advantage requires further prospective validation.

The controller leaves SCIP's source code and the
original search trajectory unchanged. Its practical
value nevertheless depends on the computing environment:
root observations and process launches incur overhead,
and savings in active worker occupancy free capacity
only when the scheduler can reallocate idle cores.
Repeated seed blocks, a prospective random gate,
and reservation-aware scheduling would help distinguish
the information value of the root statistic from
execution variability.


\section{Conclusion}
\label{sec:conclusion}

Independent-seed replication exploits MIP performance
variability but commits substantial computational resources.
We developed a finite-deadline analysis of delayed
replication and imperfect allocation decisions, examined
replication opportunities across 101,615 seeded SCIP runs,
and implemented a solver-external root-gated controller
that preserves the original search trajectory.

In a prospective two-hour evaluation on 130 paired
instances, Root-Gated proves optimality on 77 instances,
compared with 73 for both single-seed solving and
unconditional four-seed racing, while reducing active
worker occupancy by 43.00\% relative to unconditional
racing. Although these results do not isolate the causal
value of the root statistic, they demonstrate the
feasibility of selective replication under a finite
computational budget.

The broader insight is that solver variability is not
only a phenomenon to measure or mitigate, but also
an opportunity to allocate computational resources
more intelligently.

\subsubsection{\ackname} The authors used AI tools for language editing; all mathematical claims, empirical measurements, and scientific conclusions remain the sole responsibility of the authors.


\appendix
\section{Mathematical Proofs}
\label{sec:appendix_proofs}

\subsection{Proof of Proposition~\ref{prop:delayed_value}}

For a nonnegative random variable $V$, the conditional
survival-function identity gives
\begin{equation}
    \mathbb E[\min\{V,b\}\mid\mathcal I]
    =
    \int_0^b
    \Pr(V>t\mid\mathcal I)\,dt.
    \label{eq:appendix_survival}
\end{equation}

Applying this identity to $W_0$ and $W_1$ yields
\begin{equation}
\begin{aligned}
    \mathbb E[W_0\mid\mathcal I]
    &=\int_0^b S_0(t)\,dt,\\
    \mathbb E[W_1\mid\mathcal I]
    &=\int_0^b S_1(t)\,dt.
\end{aligned}
\end{equation}
Since $S_1(t)=S_0(t)$ for $t<h$,
\begin{equation}
    \mathbb E[W_1-W_0\mid\mathcal I]
    =
    -\int_h^b[S_0(t)-S_1(t)]\,dt.
    \label{eq:appendix_wall_difference}
\end{equation}

The additional replicas consume processor time for
$(W_1-h)_+$. Because $0\leq W_1\leq b$,
\begin{equation}
\begin{aligned}
    \mathbb E[(W_1-h)_+\mid\mathcal I]
    &=\int_h^b
      \Pr(W_1>t\mid\mathcal I)\,dt\\
    &=\int_h^b S_1(t)\,dt.
\end{aligned}
\label{eq:appendix_replica_work}
\end{equation}

The difference in deadline-failure indicators is
\begin{equation}
\begin{aligned}
&\mathbf 1_{\{R>b,\ h+Y>b\}}
 -\mathbf 1_{\{R>b\}}\\
&\qquad=
-\mathbf 1_{\{R>b,\ h+Y\leq b\}}.
\end{aligned}
\end{equation}
Furthermore,
\begin{equation}
\begin{aligned}
&\Pr(R>b,\ h+Y\leq b\mid\mathcal I)\\
&\quad=
\Pr(R>b\mid\mathcal I)
-\Pr(R>b,\ h+Y>b\mid\mathcal I)\\
&\quad=S_0(b)-S_1(b).
\end{aligned}
\label{eq:appendix_failure_difference}
\end{equation}

Substituting
\eqref{eq:appendix_wall_difference},
\eqref{eq:appendix_replica_work}, and
\eqref{eq:appendix_failure_difference}
into the definition of $\Delta$ in
\eqref{eq:incremental_cost}, and taking conditional
expectations, proves
\eqref{eq:conditional_replication_value}.
\hfill$\square$

\subsection{Proof of Proposition~\ref{prop:imperfect_cost}}

Because $G$ and $G^\star$ are
$\mathcal I$-measurable,
the tower property gives
\begin{equation}
\begin{aligned}
&\mathbb E[J_G-J_{G^\star}]\\
&\quad=
\mathbb E[(G-G^\star)\Delta]\\
&\quad=
\mathbb E[(G-G^\star)
          \mathbb E[\Delta\mid\mathcal I]]\\
&\quad=
\mathbb E[(G-G^\star)\delta(\mathcal I)].
\end{aligned}
\end{equation}

If $\delta(\mathcal I)<0$, then $G^\star=1$ and
\begin{equation}
    (G-G^\star)\delta(\mathcal I)
    =
    |\delta(\mathcal I)|
    \mathbf 1_{\{G\ne G^\star\}}.
\end{equation}
If $\delta(\mathcal I)\geq 0$, then $G^\star=0$
and the same identity holds. In particular, when
$\delta(\mathcal I)=0$, both sides vanish.

Thus, pointwise,
\begin{equation}
    (G-G^\star)\delta(\mathcal I)
    =
    |\delta(\mathcal I)|
    \mathbf 1_{\{G\ne G^\star\}}.
\end{equation}
Taking expectations establishes
\eqref{eq:gate_regret}.
\hfill$\square$

\subsection{Derivation of the Two-Replica Gini Identity}

For independent copies $X_1$ and $X_2$ of a
nonnegative random variable with
$0<\mu=\mathbb E[X]<\infty$,
\begin{equation}
    \min\{X_1,X_2\}
    =
    \frac{X_1+X_2-|X_1-X_2|}{2}.
\end{equation}
Taking expectations,
\begin{equation}
    \mathbb E[M_2]
    =
    \mu-\frac12\mathbb E|X_1-X_2|.
\end{equation}
By the definition
\[
    G_X=\frac{\mathbb E|X_1-X_2|}{2\mu},
\]
we obtain
\begin{equation}
    \mathbb E[M_2]=\mu(1-G_X).
\end{equation}

Under simultaneous execution, the two replicas
consume aggregate processor time $C_2^\infty=2M_2$.
Consequently,
\begin{equation}
    \mathbb E[C_2^\infty]
    =
    2\mu(1-G_X).
\end{equation}
Since $\mu>0$,
\begin{equation}
    \mathbb E[C_2^\infty]<\mu
    \quad\Longleftrightarrow\quad
    G_X>\frac12.
\end{equation}
The same derivation applies to
$X^T=\min\{X,T\}$ for a finite deadline $T$,
provided its mean is positive.
\hfill$\square$

\section{Additional Empirical and Implementation Details}
\label{sec:empirical_appendix}
\begin{table}[h]
\centering
\caption{Historical SCIP corpus and measurement availability.}
\label{tab:corpus_audit}
\begin{tabular}{lr}
\hline
Quantity & Count \\
\hline
Distinct MIP instances & 1,017 \\
Run records & 101,615 \\
Optimal termination & 42,877 \\
Time-limit termination & 56,452 \\
Other terminal statuses & 1,704 \\
Terminal status unavailable & 582 \\
First-LP basis and reduced-cost counts & 92,548 \\
Defined root degeneracy ratio & 92,048 \\
\hline
\end{tabular}
\end{table}
\begin{figure}[h]
\centering
\includegraphics[width=\linewidth]{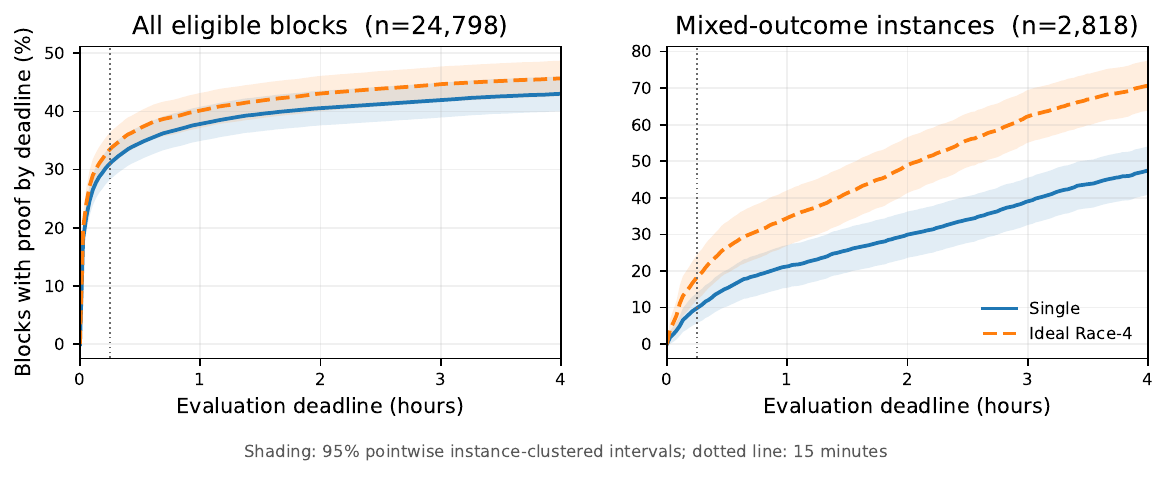}
\caption{Historical proof fractions by evaluation deadline for Single and ideal simultaneous Race-4, over all 24,798 blocks and 2,818 retrospectively mixed-outcome blocks. Shading gives 95\% pointwise instance-clustered intervals; dotted line marks 15 minutes.}
\label{fig:replication_opportunity}
\end{figure}
\begin{algorithm}[h]
\caption{Root-Gated Race-4}
\label{alg:gated_racing}
\begin{algorithmic}[1]
\Require MIP instance \(M\), seeds \(s_0,\ldots,s_3\),
deadline \(T\), threshold \(d^\star\),
polling interval \(\rho\)
\Ensure Trial summary and available worker records

\State \(t_0\gets\Call{Clock}{}\), \(D\gets t_0+T\)
\State Launch baseline process \(P_0\) with seed \(s_0\)
\State \(\mathcal P\gets\{P_0\}\),
       \(\mathrm{decided}\gets\mathrm{false}\)

\While{\(\Call{Clock}{}<D\)}
    \State Read available terminal records from
           processes in \(\mathcal P\)

    \If{a certified proof is recorded before \(D\)}
        \State Terminate remaining active processes
        \State \textbf{return} trial summary
    \EndIf

    \If{\(\neg\mathrm{decided}\) and the root record
        of \(P_0\) is available}
        \State \(\mathrm{decided}\gets\mathrm{true}\)
        \If{root record is valid and
            \(d_{\mathrm{root}}\geq d^\star\)
            and \(P_0\) remains active and unfinished}
            \For{\(k=1,2,3\)}
                \If{\(\Call{Clock}{}<D\)}
                    \State Launch \(P_k\) with seed \(s_k\)
                    \State \(\mathcal P\gets
                           \mathcal P\cup\{P_k\}\)
                \EndIf
            \EndFor
        \EndIf
    \EndIf

    \If{all processes in \(\mathcal P\) have exited}
        \State \textbf{break}
    \EndIf

    \State Poll again after interval \(\rho\)
\EndWhile

\State Wait for remaining processes to exit
       under their solver time limits
\State Record available results and resource usage
\State \textbf{return} trial summary

\end{algorithmic}
\end{algorithm}

\begin{figure}[tbp]
\centering
\includegraphics[width=\linewidth]
{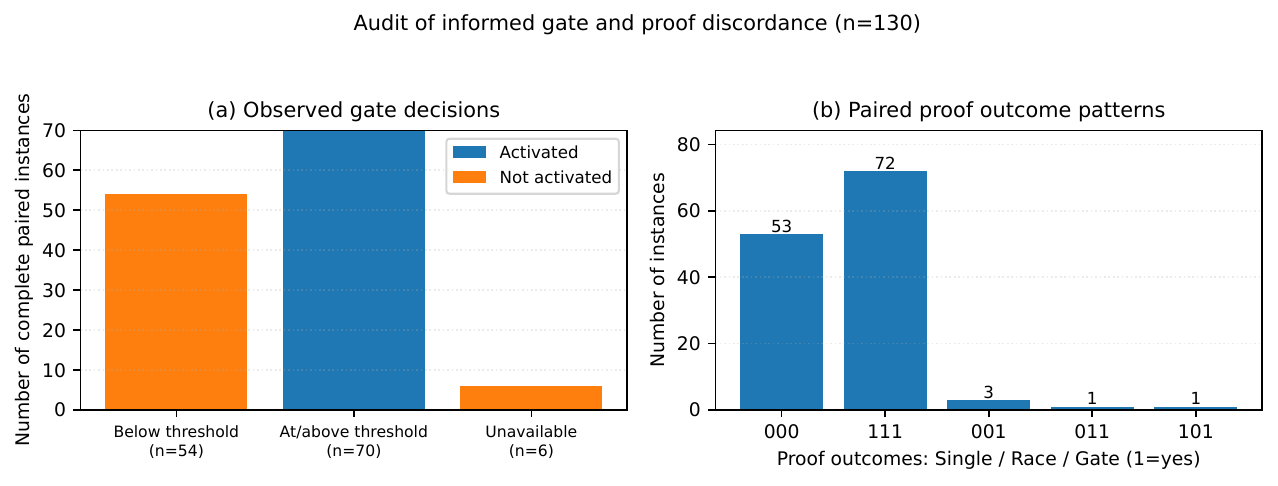}
\caption{Prospective root-gate decisions by
baseline root-LP stratum (left) and paired
optimality-proof outcome patterns for Single,
Race-4, and Root-Gated (right), across 130
complete paired instances at $T=7200$ seconds.}
\label{fig:gate_decision_audit}
\end{figure}


\begin{thebibliography}{99}
\providecommand{\natexlab}[1]{#1}
\providecommand{\url}[1]{\texttt{#1}}
\providecommand{\urlprefix}{URL }


\bibitem[{Adegbindin et~al.(2026a)}]{Adegbindin2026VariabilityScore}
Adegbindin M, Dietrich B, Pender J, Lodi A (2026) When Solvers Get Moody: Quantifying and Predicting Performance Variability in Mixed Integer Programming.
  Preprint, arXiv:2610.12876 [math.OC], \url{https://arxiv.org/abs/2610.12876}.

\bibitem[{Adegbindin et~al.(2026b)}]{Adegbindin2026LPStateSteering}
Adegbindin M, Dietrich B, Pender J (2026)
Mitigating Performance Variability in Mixed-Integer
Programming Solvers via LP State Steering.
Preprint, arXiv:2610.17268 [math.OC].
\url{https://arxiv.org/abs/2610.17268}.

\bibitem[{Anderson et~al.(2019)}]{anderson2019clairvoyant}
Anderson D, Hendel G, Le Bodic P, Viernickel M (2019)
Clairvoyant restarts in branch-and-bound search using
online tree-size estimation.
\emph{Proceedings of the AAAI Conference on Artificial
Intelligence} 33(1):1427-1434.
doi: 10.1609/aaai.v33i01.33011427.

\bibitem[{Carvajal et~al.(2014)}]{carvajal2014parallel}
Carvajal R, Ahmed S, Nemhauser G, Furman K, Goel V,
Shao Y (2014) Using diversification, communication
and parallelism to solve mixed-integer linear programs.
\emph{Operations Research Letters} 42(2):186-189.
doi: 10.1016/j.orl.2013.12.012.

\bibitem[{Danna(2008)}]{danna2008}
Danna E (2008) Performance variability in mixed integer programming.
  Presentation at the \emph{MIP 2008 Workshop on Mixed Integer Programming}, Columbia University, New York, NY.
  \url{https://coral.ise.lehigh.edu/mip-2008/talks/danna.pdf}


\bibitem[{Fischetti et~al.(2016)}]{fischetti2016random}
Fischetti M, Lodi A, Monaci M, Salvagnin D, Tramontani A (2016)
  Improving branch-and-cut performance by random sampling.
  \emph{Mathematical Programming Computation} 8(1):113-132.
  doi: 10.1007/s12532-015-0096-0.

\bibitem[{Fischetti and Monaci(2014)}]{fischetti2014erraticism}
Fischetti M, Monaci M (2014) Exploiting erraticism
in search. \emph{Operations Research} 62(1):114-122.
doi: 10.1287/opre.2013.1231.

\bibitem[{Gamrath et~al.(2020)}]{gamrath2020degeneracy}
Gamrath G, Berthold T, Salvagnin D (2020)
An exploratory computational analysis of dual
degeneracy in mixed-integer programming.
\emph{EURO Journal on Computational Optimization}
8:241-261.
doi: 10.1007/s13675-020-00130-z.

\bibitem[{Garcia(2023)}]{halSortModel}
Garcia R (2023) Performance variability in MILP modeling.
  \emph{ROADEF 2023: 24ème congrès de la Société Française de Recherche Opérationnelle et d'Aide à la Décision}, Rennes, France. \url{https://roadef2023.sciencesconf.org/434862/GARCIA_PerformanceVariability_ROADEF2023.pdf}.


\bibitem[{Gleixner et~al.(2021)}]{GleixnerEtAl2021_MIPLIB2017}
Gleixner A, Hendel G, Gamrath G, Achterberg T,
Bastubbe M, Berthold T, Christophel PM, Jarck K,
Koch T, Linderoth J, L\"ubbecke M, Mittelmann HD,
Ozyurt D, Ralphs TK, Salvagnin D, Shinano Y (2021)
MIPLIB 2017: Data-driven compilation of the 6th
mixed-integer programming library.
\emph{Mathematical Programming Computation}
13(3):443--490.
doi: 10.1007/s12532-020-00194-3.

\bibitem[{Gomes and Selman(2001)}]{gomes2001portfolios}
Gomes CP, Selman B (2001) Algorithm portfolios.
\emph{Artificial Intelligence} 126(1-2):43-62.
doi: 10.1016/S0004-3702(00)00081-3.

\bibitem[{Gurobi Optimization(n.d.)}]{gurobiConcurrentMIP}
Gurobi Optimization, LLC:
Gurobi Optimizer Reference Manual,
ConcurrentMIP parameter.
\url{https://docs.gurobi.com/projects/optimizer/en/current/reference/parameters.html}
(Accessed: 11 November 2026).



\bibitem[{Koch et~al.(2011)}]{KochEtAl2011_MIPLIB2010}
Koch T, Achterberg T, Andersen E, Bastert O, Berthold T, Bixby RE, Danna E,
Gamrath G, Gleixner AM, Heinz S, Lodi A, Mittelmann H, Ralphs T, Salvagnin D,
Steffy DE, Wolter K (2011) MIPLIB 2010: Mixed integer programming library
version 5. \emph{Mathematical Programming Computation} 3(2):103-163. doi: 10.1007/s12532-011-0025-9.


\bibitem[{Lodi and Tramontani(2014)}]{LodiTramontani2013}
Lodi A, Tramontani A (2014)
Performance variability in mixed-integer programming.
\emph{INFORMS Tutorials in Operations Research},
1-12.
doi: 10.1287/educ.2013.0112.

\end{thebibliography}
\end{document}